\documentstyle[12pt]{article} 
\newcommand{\be}{\begin{equation}}
\newcommand{\ee}{\end{equation}}
\newcommand{\bea}{\begin{eqnarray}}
\newcommand{\eea}{\end{eqnarray}}
\newcommand{\bear}{\[\begin{array}{lcl}}
\newcommand{\eear}{\end{array}\]}
\newcommand{\bearn}{\begin{equation}\begin{array}{lcl}}
\newcommand{\eearn}{\end{array}}
\newcommand{\MX}[1]{{}_{#1}M}
\newcommand{\M}[2]{{}_{#1}M_{#2}}
\newcommand{\R}[1]{R_{{}_{[#1]}}}
\newcommand{\RX}[1]{{}_{#1}R}
\newcommand{\ST}[3]{S_{{}_{(#1)(#2)(#3)}}}

\newcommand{\Q}[2]{Q_{{}_{(#1)(#2)}}}

\def\emline#1#2#3#4#5#6{%
       \put(#1,#2){\special{em:moveto}}%
       \put(#4,#5){\special{em:lineto}}}

\begin{document} 

\title{A note on the eight tetrahedron equations}
\author{Jarmo Hietarinta\thanks{E-mail: hietarin@newton.tfy.utu.fi}\\
Department of Physics, University of Turku\\ FIN-20014 Turku, Finland\\
and\\
Frank Nijhoff\thanks{E-mail: frank@amsta.leeds.ac.uk}\\
Department of Applied Mathematical Studies\\ 
University of Leeds, Leeds LS2 9JT,UK} 

\maketitle

\begin{abstract}
In this paper we derive from arguments of string scattering
a set of eight tetrahedron equations, with different
index orderings. It is argued that this system of equations is the
proper system that represents integrable structures in three
dimensions generalising the Yang-Baxter equation. Under additional
restrictions this system reduces to the usual tetrahedron equation in
the vertex form. Most known solutions fall under this class, but it is
by no means necessary. Comparison is made with the work on braided
monoidal 2-categories also leading to eight tetrahedron equations.
\end{abstract}

\vfill\eject

\section{Introduction}
An important current problem in the study of integrable systems is to
make an extension to higher dimensions.  For 1+1 dimensions there are
several well established approaches to integrability and many
beautiful results have been obtained; much less is known about
integrability of three dimensional systems.  One of the most important
approaches to 1+1 dimensional integrability is based on the
Yang-Baxter equations \cite{YBE}, the corresponding 2+1 dimensional
``tetrahedron'' equations were introduced by Zamolodchikov already in
the 1980 \cite{Zam80,Zam81,BS}.  These equations have been under
intense study during the last few years
\cite{MN,MN1,MN2,Baz1,Baz2,Kor1,Kor2,Labels,SMS,Hor}, but many
fundamental questions still remain open.

One difference between 1+1 and 2+1 dimensional integrability stems
from the fact that there is no natural ordering in the two dimensional
space.  The Yang-Baxter equation can be derived, e.g., from the
condition of factorizable scattering of point particles on a line
\cite{YBE} and since one can introduce a good ordering on a line there
is no ambiguity in writing down the Yang-Baxter equations.
Zamolodchikov's tetrahedron equations can be derived from the
conditions of factorizable scattering of straight strings
\cite{Zam80,Labels} (or particles at the intersections of strings
\cite{Labels}) on a plane.  In the particle interpretation the
scattering matrix depends on three incoming and three outgoing
particles, but since there is no obvious way of defining an order in
two dimensions, there is no single ordering in which the indices of
the corresponding scattering matrix should be written.

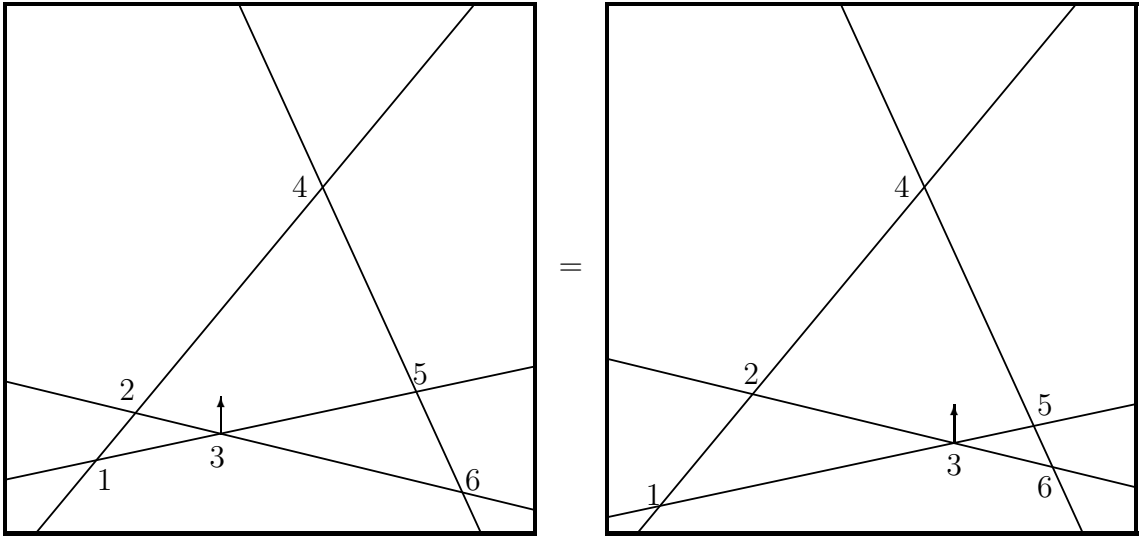
\begin{figure}[t]
\unitlength=1mm
\special{em:linewidth 0.7pt}
\linethickness{1.2pt}
\begin{picture}(150,70)(-5,0)
\put(80,0){\framebox(70,70)[cc]{}}
\put(0,0){\framebox(70,70)[cc]{}}
\linethickness{0.6pt}
\emline{84}{0}{1}{142}{70}{2}\emline{143}{0}{3}{111}{70}{4}
\emline{4}{0}{5}{62}{70}{6}\emline{0}{7}{7}{70}{22}{8}
\emline{70}{3}{9}{0}{20}{10}\emline{63}{0}{11}{31}{70}{12}
\emline{80}{2}{13}{150}{17}{14}\emline{80}{23}{15}{150}{6}{16}
\put(28.5,13){\vector(0,1){5}}\put(126,12){\vector(0,1){5}}
\put(75,35){\makebox(0,0)[cc]{=}}
\put(13,7){\makebox(0,0)[cc]{1}}\put(16,19){\makebox(0,0)[cc]{2}}
\put(39,46){\makebox(0,0)[cc]{4}}\put(28,10){\makebox(0,0)[cc]{3}}
\put(55,21){\makebox(0,0)[cc]{5}}\put(62,7){\makebox(0,0)[cc]{6}}
\put(86,5){\makebox(0,0)[cc]{1}}\put(99,21){\makebox(0,0)[cc]{2}}
\put(126,9){\makebox(0,0)[cc]{3}}\put(119,46){\makebox(0,0)[cc]{4}}
\put(138,17){\makebox(0,0)[cc]{5}}\put(138,6){\makebox(0,0)[cc]{6}}
\end{picture}
\caption{The starting configuration for four string (six particle) 
scattering. The resulting total scattering matrix should be the same
for the two alternatives that differ only in the relative position of
the two strings at the bottom.}
\label{F:1}
\end{figure}

Let first recall how the tetrahedron equation arises when we consider
the scattering of straight strings. The basic scattering process is
that of three straight strings, and if we are only interested in the
particles at the intersection of the strings (the ``vertex''
formulation) the scattering matrix is written as
$S_{j_1j_2j_3}^{k_1k_2k_3}$ where the $j$'s give the states of
incoming particles and the $k$' the states of the outgoing ones. The
tetrahedron equations arise when we consider the scattering of four
strings \cite{Zam80,Labels}, which generically have six intersections,
see Figure \ref{F:1}.
The initial configuration looks like an arrow, and if we go to frame
where the arrowhead (particle 4) is stationary, the dynamics is
described fully by the way the intersection point 3 moves. Depending
on the relative initial positions of the two strings at the bottom of
the figure, particle 3 will pass particle 4 on the left or on the
right. In each alternative there will be four basic scattering
processes, in each of which a triangle will be turned over.  These two
alternatives should give the same result, and this condition yields
the tetrahedron equations:
\be
S_{j_1j_2j_3}^{k_1k_2k_3}S_{k_1j_4j_5}^{l_1k_4k_5}
S_{k_2k_4j_6}^{l_2l_4k_6}S_{k_3k_5k_6}^{l_3l_5l_6}=
S_{j_3j_5j_6}^{k_3k_5k_6}S_{j_2j_4k_6}^{k_2k_4l_6}
S_{j_1k_4k_5}^{k_2l_4l_5}S_{k_1k_2k_3}^{l_1l_2l_3}.
\ee
Here we have used Einstein summation convention for the repeated $k$
indices.

In writing down the above equation we have used a particular
convention for the index ordering: for each triangle that turns over
we have taken the indices from left to right. Since the four string
configuration of Fig.\ \ref{F:1} is not rotationally invariant
``left'' can always be defined, but any such ordering gives problems
already when one considers the scattering of five strings.

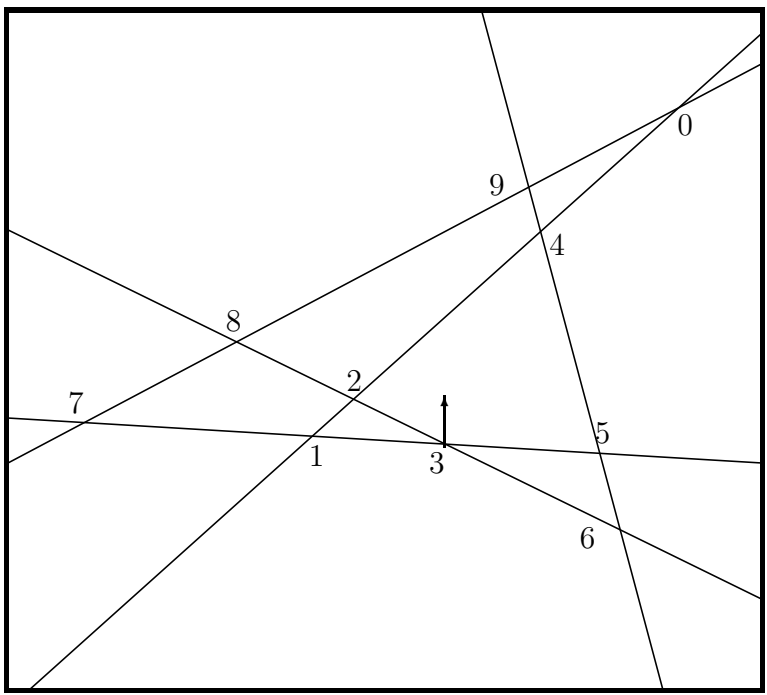
\begin{figure}[t]
\unitlength=1mm
\special{em:linewidth 0.6pt}
\linethickness{1.5pt}
\begin{picture}(120,100)(-9,30)
\put(20,30){\framebox(100,90)[cc]{}}
\linethickness{0.6pt}
\emline{120}{60}{1}{20}{66}{2}
\put(78,62){\vector(0,1){7}}
\emline{20}{60}{3}{120}{113}{4}\emline{120}{42}{5}{20}{91}{6}
\emline{23}{30}{7}{120}{117}{8}\emline{107}{30}{9}{83}{120}{10}
\put(61,61){\makebox(0,0)[cc]{1}}\put(66,71){\makebox(0,0)[cc]{2}}
\put(77,60){\makebox(0,0)[cc]{3}}\put(93,89){\makebox(0,0)[cc]{4}}
\put(99,64){\makebox(0,0)[cc]{5}}\put(97,50){\makebox(0,0)[cc]{6}}
\put(29,68){\makebox(0,0)[cc]{7}}\put(50,79){\makebox(0,0)[cc]{8}}
\put(85,97){\makebox(0,0)[cc]{9}}\put(110,105){\makebox(0,0)[cc]{0}}
\end{picture}
\caption{A typical situation with five string scattering. The first 
triangle to turn over is 123 of the arrow 1463, but if it is
considered as part of the arrow 3702 then the corresponding scattering
matrix should be labeled as 312.}
\label{F:2}
\end{figure}

A typical starting configuration of a five string scattering is given
in Figure \ref{F:2}.
Let us assume that in the next scattering process the triangle 123
will be turned over. In which order should we now write the indices of
the corresponding scattering matrix? If we consider the triangle 123
as a part of arrow 1463 we should use $S_{123}$, according to the
above convention, but if we consider arrow 3702 and look at the
picture from right, we should use $S_{312}$.  This problem was
recognized in \cite{Zam80} and was taken care of by requiring that the
$S$-matrix is invariant under cyclic index permutation,
see. Eq. (3.5) of \cite{Zam80}.

In this paper we show that this ordering ambiguity means that there
are, in fact, eight tetrahedron equations (obtained from the standard
one by certain index permutations) that must be satisfied
simultaneously by the tetrahedron $S$-matrix. In Sec.\ 2 we give an
algebraic derivation of these equations using the ``obstruction''
method, cf. \cite{MN1,MN2,Kor1,Hor}. In Sec.\ 3 we give several
interpretations to these equations and discuss the conditions under
which these equations collapse into one.  In Sec.\ 4 we will make a
connection with the notion of higher Bruhat orders introduced by Manin
and Schechtman in \cite{MS,MS2}. Another formulation of the tetrahedron
equations is in terms of braided monoidal bicategories,
cf. \cite{KV1,KV2}, and provides an alternative way of obtaining the
system of eight tetrahedron equations \cite{VK}. However, we believe
our derivation is closer to the physical interpretation in terms of
string scattering, furthermore we will not need to use the language of
bicategories. Since our derivation relies on the obstruction
mechanism, \cite{MN1,MN2}, we hope that eventually this point of view
leads to the derivation (in the spirit of \cite{Kor1,Kor2}) of
explicit solutions of the system in the cases when it does not
collapse to a single tetrahedron equation.

\section{Derivation} 
\subsection{Derivation of the Yang-Baxter equation}
\label{S:DYB}
We will start by recalling the algebraic derivation of the YBE.  Let
us assume that we have a set of $d$ $n\times n$ matrices which also
depend on a continuous `spectral' (or `color') parameters: ${\cal
M}=\{\MX i(\lambda ) \in End(V_{\emptyset})| i=1,\dots,d,\,\lambda\in
\bf CP\}$, using the convention that the matrix indices are written on
the right and the other indices on the left. For later purposes it is
useful to think of the spectral parameter $\lambda$ as being some
projective vector over $\bf C$. Let us now assume that the matrices of
$\cal M$ do not quite commute but that their commutation is
``obstructed'' by some numerical coefficients $R$: 
\be
\MX{j_{1}}(\lambda_{1})_{\alpha}^{\beta}\;
\MX{j_{2}}(\lambda_{2})_{\beta}^{\gamma}=
R(\lambda_{1},\lambda_{2})_{j_{1}j_{2}}^{k_{1}k_{2}}\;
\MX{k_{2}}(\lambda_{2})_{\alpha}^{\beta}\;
\MX{k_{1}}(\lambda_{1})_{\beta}^{\gamma}.
\label{E:YBr1}
\ee
The obstruction coefficients $R$ can be put into a $d^2\times d^2$
matrix and we can say that it operates on the product of vector spaces
$V_1\otimes V_2$, whose basis is given by the matrices $\MX i$,
themselves operating on some other vector space $V_{\emptyset}$. [This
hierarchical structure will be taken one step further in the next
section.] We can now use a shorthand notation and write down only the
names of the vector spaces where the operation takes place
\be
\M{1}\emptyset\cdot\M{2}\emptyset = 
R_{12}\;\M{2}\emptyset\cdot\M{1}\emptyset.
\label{E:YBr2}
\ee
It should be remembered that {\em with each vector space comes its own
spectral parameter} (the parameter associated with $V_{\emptyset}$ is
global).

[There is an alternative way of obstructing commutativity by
\[
R_{ik}^{pq}\;(\lambda,\mu)\;{}_p^mT(\lambda)\cdot{}_q^nT(\mu)=
{}_k^rT(\mu)\cdot{}_i^sT(\lambda)\;R_{rs}^{mn}(\mu,\lambda),
\]
where the ${}_i^jT$'s are some non-commuting quantities, each of which
can be represented by a matrix acting on some vector space.
Multiplying by $R^{-1}$ from the left we can write this in the form
(\ref{E:YBr1}), but with double indices.]

If the reversal (\ref{E:YBr1}) is done twice we get
\[
\left[\delta_1\delta_2 - R_{12}R_{21}\right]\M{1}\emptyset
\cdot\M{2}\emptyset=0,
\]
which is usually taken in the strong form as the `unitarity' condition
\be
R_{12}R_{21}=\delta_1\delta_2.
\label{E:YBuni}
\ee

Taking into account the associativity of the matrix product we see that 
the obstruction to commutativity (\ref{E:YBr2}) leads to two different 
ways of inverting the triple $ABC$, namely on the one hand: 
$(AB)C \to B(AC) \to (BC)A \to CBA$, and on the other hand 
$A(BC) \to (AC)B \to C(AB) \to CBA$. Equating the two expressions 
obtained by elaborating these two ways, namely 
\bea
\M1\emptyset\cdot\M2\emptyset\cdot\M3\emptyset&=&R_{12}\;
\M2\emptyset\cdot\M1\emptyset\cdot\M3\emptyset\nonumber\\
&=&R_{12}R_{13}\;\M2\emptyset\cdot\M3\emptyset\cdot\M1\emptyset\nonumber\\
&=&R_{12}R_{13}R_{23}\;\M3\emptyset\cdot\M2\emptyset
\cdot\M1\emptyset,\nonumber
\eea
and 
\bea
\M1\emptyset\cdot\M2\emptyset\cdot\M3\emptyset&=&R_{23}\;
\M1\emptyset\cdot\M3\emptyset\cdot\M2\emptyset\nonumber\\
&=&R_{23}R_{13}\;\M3\emptyset\cdot\M1\emptyset\cdot\M2\emptyset\nonumber\\
&=&R_{23}R_{13}R_{12}\;\M3\emptyset\cdot\M2\emptyset\cdot
\M1\emptyset.\nonumber
\eea
we obtain in the strong sense the quantum Yang-Baxter equation as a
condition on $R$:
\be
R_{12}R_{13}R_{23}=R_{23}R_{13}R_{12},
\label{E:YBsh}
\ee
which is short-hand for 
\be
R(\lambda_1,\lambda_2)_{j_1j_2}^{k_1k_2}
R(\lambda_1,\lambda_3)_{k_1j_3}^{l_1k_3}
R(\lambda_2,\lambda_3)_{k_2k_3}^{l_2l_3}=
R(\lambda_2,\lambda_3)_{j_2j_3}^{k_2k_3}
R(\lambda_1,\lambda_3)_{j_1k_3}^{k_1l_3}
R(\lambda_1,\lambda_2)_{k_1k_2}^{l_1l_2}\   .
\label{E:YBlo}
\ee

\subsection{Derivation of the tetrahedron equation}
We will next derive in a similar way the tetrahedron equation
\cite{BS,MN}. We start with an indexed set of matrices $R$ operating on
a product of two vector spaces, i.e.\ ${\cal R}= \{\RX{\alpha_{ij}}_{ij} 
(\lambda_i,\lambda_j)\in End(V_i,V_j)|i,j=1,\dots,n,\; \alpha_{ij}=1,
\dots,m,\, \lambda_i\in {\bf CP^2}\}$, where the spectral parameter is 
now a projective 3-dimensional vector. As an extension of the previous 
case, we assume that the $R$'s do not quite satisfy the Yang-Baxter
equation, but rather obey 
\bea
&&\RX{\alpha_{12}}_{12}(\lambda_1,\lambda_2)\cdot
\RX{\alpha_{13}}_{13}(\lambda_1,\lambda_3)\cdot
\RX{\alpha_{23}}_{23}(\lambda_2,\lambda_3)=\nonumber\\
&&\quad S(\lambda_1,\lambda_2,\lambda_3)_{\alpha_{12}
\alpha_{13}\alpha_{23}}^{\beta_{12}\beta_{13}\beta_{23}}\;
\RX{\beta_{23}}_{23}(\lambda_2,\lambda_3)\cdot
\RX{\beta_{13}}_{13}(\lambda_1,\lambda_3)\cdot
\RX{\beta_{12}}_{12}(\lambda_1,\lambda_2),
\label{E:Zobs}
\eea
which defines the obstruction matrix $S$, operating on the product of
three vector spaces $V_{(12)}\otimes V_{(13)}\otimes V_{(23)}$,
labeled now by {\it pairs} of integers.  In (\ref{E:Zobs}) the internal
indices of the $R$'s have been indicated only by the vector spaces on
which they act, and there is a distributed matrix product just like in the
Yang-Baxter equation over them. The external indices $\alpha_{ij},\,
\beta_{ij}$ are written out explicitly, and there is a summation over
the $\beta_{ij}$'s.

Korepanov has successfully used (\ref{E:Zobs}) in constructing
solutions to the tetrahedron equation \cite{Kor1,Kor2}, by choosing
suitably deformed solutions of the Yang-Baxter equation, see also
\cite{Hor}.

Since the left and right hand sides of (\ref{E:YBlo},\ref{E:Zobs})
have different index distributions (observe the positions of the
repeated indices in the $R$'s) we will also need another
reversal
\bea
&&\RX{\alpha_{23}}_{23}(\lambda_2,\lambda_3)\cdot
\RX{\alpha_{13}}_{13}(\lambda_1,\lambda_3)\cdot
\RX{\alpha_{12}}_{12}(\lambda_1,\lambda_2)=\nonumber\\
&&\quad \tilde S(\lambda_3,\lambda_2,\lambda_1)_{\alpha_{23}
\alpha_{13}\alpha_{12}}^{\beta_{23}\beta_{13}\beta_{12}}\;
\RX{\beta_{12}}_{12}(\lambda_1,\lambda_2)\cdot
\RX{\beta_{13}}_{13}(\lambda_1,\lambda_3)\cdot
\RX{\beta_{23}}_{23}(\lambda_2,\lambda_3).
\label{E:Zobs2}
\eea
To simplify notation we will only write down the indices of
the various spaces
\be
\R{12}\cdot \R{13}\cdot \R{23}
=\ST{12}{13}{23}\;\R{23}\cdot \R{13}\cdot \R{12},
\label{E:rev1}
\ee
\be
\R{23}\cdot \R{13}\cdot \R{12}
=\tilde \ST{23}{13}{12}\;\R{12}\cdot \R{13}\cdot \R{23}.
\label{E:rev2}
\ee
Here the square brackets around the subscripts of $R$ are to remind us
that there are both external indices labeling the different $R$
matrices and internal indices of Yang-Baxter type, while the brackets
around the indices of $S$ and $\tilde S$ indicate that they are only
external indices. Note that the order inside each bracketed pair is
also important, and relabelings should be made with caution.

So far we have no relation between $S$ and $\tilde S$, because they
arose from different reversals. However, an application of these two
reversals in succession yields the starting order, suggesting that the
unitarity condition
\be
\tilde \ST{23}{13}{12}\ST{12}{13}{23}=
\delta_{{}_{(12)}}\delta_{{}_{(13)}}\delta_{{}_{(23)}},
\label{E:Uni}
\ee
should be satisfied, but again this is necessary only in the weak
sense, i.e. when acting on a triple of matrices $R$.

In addition to the above we have to give a rule for exchanging $R$'s
with disjoint indices. In general we could introduce a permutation
operator $Q$ by \cite{Ruth}
\be
\R{12}\; \R{34}=\Q{12}{34}\;\R{34}\;\R{12}.
\label{E:RLE1}
\ee Since the internal (matrix) indices are disjoint it would be
natural to choose $Q=\delta\delta$, but even if we later may take this
conventional choice it is useful to carry along the operator $Q$,
since it will turn out to be a good place-marker in the otherwise
monotonous tetrahedron equation.  Furthermore, Lawrence has proposed
in \cite{Ruth} a variant of the tetrahedron equation where this $Q$
operator is taken into account. [Among other things this allows for
some additional (reductive) solutions of the form: ~$Q_{(12)(34)}$~ a
solution of the quantum Yang-Baxter equation, while
~$S_{(12)(13)(23)}=Q_{(12)(13)} Q_{(12)(23)}Q_{(13)(23)}$~.]
In this paper we will only use the commutation and inversion 
properties:
\be 
\Q{12}{34}\Q{13}{24}=\Q{13}{24}\Q{12}{34},
\label{E:Qcomm}
\ee
\be
\Q{12}{34}\Q{34}{12}=\delta_{(12)}\delta_{(34)},
\label{E:QQ}
\ee

$S$ and $Q$ we defined by reversals of three and two $R$'s,
respectively.  When we considers ways of reversing more than three
objects we get conditions for $S$ (and $Q$). In fact, because of the
dependence on {\em pairs} of indices we need to consider next the
reversal of a product of six objects: $\R{ij}$ where $i\not=
j\in\{1,2,3,4\}$.  One particular case is the following: 
\bea &
\underbrace{\R{12}\cdot\R{13}\cdot\R{23}}
\cdot\R{14}\cdot\R{24}\cdot\R{34}&\nonumber\\ &\downarrow&
\ST{12}{13}{23}\nonumber\\ &\R{23}\cdot\R{13}\cdot\underbrace{\R{12}
\cdot\R{14}\cdot\R{24}}\cdot\R{34}&\nonumber\\ &\downarrow
&\ST{12}{14}{24}\nonumber\\ &\R{23}\cdot\underline{\R{13}\cdot\R{24}}
\cdot\R{14}\cdot\underline{\R{12}\cdot\R{34}} &\nonumber\\
&\Downarrow&\Q{13}{24}\Q{12}{34}\nonumber\\
&\R{23}\cdot\R{24}\cdot\underbrace{\R{13}
\cdot\R{14}\cdot\R{34}}\cdot\R{12}&\nonumber\\ &\downarrow&
\ST{13}{14}{34}\nonumber\\ &\underbrace{\R{23}\cdot\R{24}\cdot\R{34}}
\cdot\R{14}\cdot\R{13}\cdot\R{12}&\nonumber\\ &\downarrow
&\ST{23}{24}{34}\nonumber\\ &\R{34}\cdot\R{24}\cdot\underline{\R{23}
\cdot\R{14}}\cdot\R{13}\cdot\R{12}&\nonumber\\
&\Downarrow&\Q{23}{14}\nonumber\\
&\R{34}\cdot\R{24}\cdot\R{14}\cdot\R{23}
\cdot\R{13}\cdot\R{12}&\nonumber \eea Here the under-brace indicates
which triple is reversed by the $S$ and the underline means the terms
are commuted using $Q$. At each step the multiplying obstruction
matrix is written at the Dow-narrow, and they compose as
\[
\ST{12}{13}{23}\ST{12}{14}{24}\Q{13}{24}\Q{12}{34}
\ST{13}{14}{34}\ST{23}{24}{34}\Q{23}{14}.
\]

There is precisely one other way to reverse the previous starting point:
\bea
&\R{12}\cdot\R{13}\cdot
\underline{\R{23}\cdot\R{14}}
\cdot\R{24}\cdot\R{34}&\nonumber\\
&\Downarrow&\Q{23}{14}\nonumber\\
&\R{12}\cdot\R{13}\cdot\R{14}\underbrace{
\cdot\R{23}\cdot\R{24}\cdot\R{34}}&\nonumber\\
&\downarrow& \ST{23}{24}{34}\nonumber\\
&\R{12}\cdot\underbrace{\R{13}\cdot\R{14}
\cdot\R{34}}\cdot\R{24}\cdot\R{23}&\nonumber\\
&\downarrow& \ST{13}{14}{34}\nonumber\\
&\underline{\R{12}\cdot\R{34}}\cdot\R{14}
\cdot\underline{\R{13}\cdot\R{24}}\cdot\R{23}&\nonumber\\
&\Downarrow&\Q{12}{34}\Q{13}{24}\nonumber\\
&\R{34}\cdot\underbrace{\R{12}\cdot\R{14}
\cdot\R{24}}\cdot\R{13}\cdot\R{23}&\nonumber\\
&\downarrow& \ST{12}{14}{24}\nonumber\\
&\R{34}\cdot\R{24}\cdot\R{14}
\cdot\underbrace{\R{12}\cdot\R{13}\cdot\R{23}}
&\nonumber\\&\downarrow &\ST{12}{13}{23}\nonumber\\
&\R{34}\cdot\R{24}\cdot\R{14}
\cdot\R{23}\cdot\R{13}\cdot\R{12}&\nonumber
\eea
and for the last line the multiplier will be
\[
\Q{23}{14}\ST{23}{24}{34}\ST{13}{14}{34}\Q{12}{34}\Q{13}{24}
\ST{12}{14}{24}\ST{12}{13}{23}.
\]
The equality of the above two expressions yields the tetrahedron equation:
\bea
&&\ST{12}{13}{23}\ST{12}{14}{24}\Q{13}{24}\Q{12}{34}
\ST{13}{14}{34}\ST{23}{24}{34}\Q{23}{14}=\nonumber\\
&&\quad\quad\Q{23}{14}\ST{23}{24}{34}\ST{13}{14}{34}\Q{12}{34}\Q{13}{24}
\ST{12}{14}{24}\ST{12}{13}{23}.
\label{E:TEsl}
\eea
If $Q=\delta\delta$ and we use the translation table 1=12, 2=13, 3=23,
4=14, 5=24, 6=34, we get the tetrahedron equation in the usual
notation
\be\label{E:SSSS} 
S_{123}S_{145}S_{246}S_{356}=S_{356}S_{246}S_{145}S_{123}.
\ee
We note that the double index notation of (\ref{E:TEsl}) is more
natural, because it identifies the points by the intersections of the
two strings. 

Let us finish this section by a comment on the spectral parameters.
Each matrix $S$ depends on three spectral parameters attached to each
of the labels it carries. The derivation of this section gives a
natural distribution of these parameters through the equation, as it
does in the Yang-Baxter case. However, since we take the spectral
parameter to be a projective vector in a fixed 3-dimensional complex
space there is nonetheless a condition arising from the fact that four
vectors in a three-dimensional space are necessarily dependent.  This
leads to the determinant condition given in \cite{MN} which, when
expressed in terms of spherical angles associated with each of these
vectors, is exactly the condition between the spectral parameters used
in Zamolodchikov's construction of his solution, cf. \cite{Zam81}.  An
interesting question is whether this parametrization would correspond
to the one that one would expect from the Baxterization procedure via
the generalization of Coxeter groups underlying the tetrahedron
equations, as proposed in \cite{Bel,Viallet,Mai}

\subsection{The other tetrahedron equations}
The main observation of this paper is that the above picture is
incomplete in view of the fact that there are other starting points
for the reversal of six $R$-matrices which will lead to tetrahedron
equations which in general are not equivalent to (\ref{E:TEsl}) or
(\ref{E:SSSS}). In fact, we should investigate all possible starting
configurations of matrices $R$ and thus obtain a set of equations
involving $S$ as well as $\tilde{S}$.

It is not hard to find those starting configurations for which at
least two triple reversals can be done.  Without any loss of
generality we may renumber the vector spaces and indices so that the
first reversal is on $\cdots\R{12} \cdot \R{13} \cdot \R{23}\cdots$
resulting with $\cdots\R{23}\cdot\R{13}\cdot\R{12}\cdots$ or in the
reverse order: $\cdots\R{23}\cdot\R{13}\cdot\R{12}\cdots$ resulting
with $\cdots\R{12} \cdot \R{13} \cdot \R{23}\cdots$.  The next
reversal must involve $\R{12}$, $\R{23}$ or $\R{13}$.  In the first
case the two other $R$'s that go with $\R{12}$ must be on its right
hand side, and can be numbered as $\R{14}\cdot\R{24}$ or $\R{42}
\cdot\R{41}$ (note the order of indices,
c.f. (\ref{E:rev1},\ref{E:rev2})), and the remaining term $\R{\{34\}}$
(we do not yet know which index order works, this is reminded by the
braces) can then be put in three different places resulting with six
starting configurations: 
\bear
1\quad&&\R{12}\cdot\R{13}\cdot\R{23}\cdot
\R{14}\cdot\R{24}\cdot\R{\{34\}}\\
1'\quad&&\R{12}\cdot\R{13}\cdot\R{23}
\cdot\R{42}\cdot\R{41}\cdot\R{\{34\}}\\
2'\quad&&\R{12}\cdot\R{13}\cdot\R{23}\cdot
\R{\{34\}}\cdot\R{14}\cdot\R{24}\\
2\quad&&\R{12}\cdot\R{13}\cdot\R{23}
\cdot\R{\{34\}}\cdot\R{42}\cdot\R{41}\\
3\quad&&\R{\{34\}}\cdot\R{12}\cdot\R{13}\cdot\R{23}\cdot
\R{14}\cdot\R{24}\\ 3'\quad&&\R{\{34\}}\cdot\R{12}\cdot\R{13}\cdot
\R{23}\cdot\R{42}\cdot\R{41} 
\eear 
If the second reversal involves $\R{23}$, its (left hand side)
companions can be numbered as $\R{42}\cdot\R{43}$ or
$\R{34}\cdot\R{24}$, and the remaining term $\R{\{14\}}$ can again be
distributed among the terms in three ways. This results with the
following six possible starting configurations:
\bear 
4\quad&&\R{42}\cdot\R{43}\cdot\R{12}
\cdot\R{13}\cdot\R{23}\cdot\R{\{14\}}\\
4'\quad&&\R{34}\cdot\R{24}\cdot\R{12}
\cdot\R{13}\cdot\R{23}\cdot\R{\{14\}}\\
5'\quad&&\R{42}\cdot\R{43}\cdot\R{\{14\}}\cdot\R{12}
\cdot\R{13}\cdot\R{23}\\
5\quad&&\R{34}\cdot\R{24}\cdot\R{\{14\}}\cdot\R{12}
\cdot\R{13}\cdot\R{23}\\
6\quad&&\R{\{14\}}\cdot\R{42}\cdot\R{43}\cdot\R{12}
\cdot\R{13}\cdot\R{23}\\
6'\quad&&\R{\{14\}}\cdot\R{34}\cdot\R{24}\cdot\R{12}
\cdot\R{13}\cdot\R{23} 
\eear 
Finally, if the second reversal uses $\R{13}$ we must put $\R{14}$ on
its left hand side and $\R{43}$ of the right and $\R{\{24\}}$ on
either end, so that we can start with 
\bear 
7\quad&&\R{\{24\}}\cdot
\R{14}\cdot \R{12}\cdot \R{13}\cdot \R{23} \cdot \R{43}\\
8\quad&&\R{14}\cdot \R{12}\cdot \R{13}\cdot \R{23}\cdot \R{43}\cdot
\R{\{24\}} 
\eear

These starting points are then guaranteed to allow at least two triple
reversals, but for a complete order reversal we have to do four triple
reversals. It turns out that a third reversal cannot be done in all
cases, these bad cases are marked above with a prime.  For each of the
remaining eight starting points the first reversal can be done in
precisely two ways, one of them is by $S_{(12)(13)(23)}$, the other
one varying from case to case, these two alternatives give the two
sides of the tetrahedron equations. We will not give the details here,
the derivation follows the one done before and is easy to do since at
each step there are no alternatives in applying the triple reversals.

The above classification can be repeated for starting points for which
the first reversal is on $\cdots\R{23}\cdot\R{13}\cdot\R{12}\cdots$
resulting with $\cdots\R{12} \cdot \R{13} \cdot \R{23}\cdots$.
However, it turns out that these reversed starting points can be
relabeled so that they give the same as the unreversed ones, except
for cases 1 and 6. Thus we obtain, finally, eight different
tetrahedron equations:
\bearn
1,6\quad&&\ST{12}{13}{23}\ST{12}{14}{24}\Q{13}{24}\Q{12}{34}
\ST{13}{14}{34}\ST{23}{24}{34}\Q{23}{14}
=\\&&\quad\quad\Q{23}{14}\ST{23}{24}{34}
\ST{13}{14}{34}\Q{12}{34}\Q{13}{24}
\ST{12}{14}{24}\ST{12}{13}{23},\\
2,8r\quad&&\ST{12}{13}{23}\Q{12}{43}\tilde \ST{12}{42}{41}
\tilde \ST{13}{43}{41}\Q{23}{41}\Q{13}{42}\tilde \ST{23}{43}{42}
=\\&&\quad\quad\tilde \ST{23}{43}{42}\Q{13}{42}\Q{23}{41}
\tilde \ST{13}{43}{41}\tilde \ST{12}{42}{41}\Q{12}{43}\ST{12}{13}{23},\\
3,4r\quad&&\ST{12}{13}{23}\ST{12}{14}{24}\Q{13}{24}
\tilde \ST{43}{23}{24}\tilde \ST{43}{13}{14}\Q{23}{14}\Q{43}{12}
=\\&&\quad\quad\Q{43}{12}\Q{23}{14}\tilde \ST{43}{13}{14}\tilde 
\ST{43}{23}{24}\Q{13}{24}\ST{12}{14}{24}\ST{12}{13}{23},\\
4,3r\quad&&\ST{12}{13}{23}\ST{42}{43}{23}\Q{42}{13}
\tilde \ST{42}{12}{14}\tilde \ST{43}{13}{14}\Q{23}{14}\Q{43}{12}
=\\&&\quad\quad\Q{43}{12}\Q{23}{14}\tilde \ST{43}{13}{14}\tilde 
\ST{42}{12}{14}\Q{42}{13}\ST{42}{43}{23}\ST{12}{13}{23},\\
5,7r\quad&&\ST{12}{13}{23}\Q{14}{23}
\tilde \ST{34}{24}{23}\tilde \ST{34}{14}{13}\Q{24}{13}\Q{34}{12}
\tilde \ST{24}{14}{12}
=\\&&\quad\quad\tilde \ST{24}{14}{12}\Q{34}{12}\Q{24}{13}
\tilde \ST{34}{14}{13}\tilde \ST{34}{24}{23}\Q{14}{23}\ST{12}{13}{23},\\
1r,6r\quad&&
\tilde \ST{23}{13}{12}\tilde \ST{24}{14}{12}\Q{24}{13}\Q{34}{12}
\tilde \ST{34}{14}{13}\tilde \ST{34}{24}{23}\Q{14}{23}=\\&&\quad\quad
\Q{14}{23}\tilde \ST{34}{24}{23}\tilde \ST{34}{14}{13}
\Q{34}{12}\Q{24}{13}\tilde \ST{24}{14}{12}\tilde \ST{23}{13}{12},\\
7,5r\quad&&\ST{12}{13}{23}\Q{14}{23}\Q{12}{43}
\ST{14}{13}{43}\ST{24}{23}{43}\Q{24}{13}\tilde \ST{24}{14}{12}=
\\&&\quad\quad\tilde \ST{24}{14}{12}\Q{24}{13}
\ST{24}{23}{43}\ST{14}{13}{43}\Q{12}{43}\Q{14}{23}\ST{12}{13}{23},\\
8,2r\quad&&\ST{12}{13}{23}\Q{14}{23}\Q{12}{43}
\ST{14}{13}{43}\ST{14}{12}{42}\Q{13}{42}\tilde \ST{23}{43}{42}=
\\&&\quad\quad\tilde \ST{23}{43}{42}\Q{13}{42}
\ST{14}{12}{42}\ST{14}{13}{43}\Q{12}{43}\Q{14}{23}\ST{12}{13}{23}.
\eearn\label{E:l1}\ee

\section{Interpretation}
To get a better understanding of the equations (\ref{E:l1}), let us 
look at the  simplified case where the matrices $Q$ are all taken equal 
to one, and investigate the geometric meaning of the set of equations 
we have obtained. 

\subsection{Reduction under unitarity}
Let us first renumber the indices in (\ref{E:l1}) so that inside each
bracket $(ij)$ we have $i<j$.  This is accomplished by the following
cyclic renumberings: 1 none, 2 (1234), 3 (34), 4 (234), 5 none, 1r
none, 7 (34), 8 (234).  After this it turns out that in each $S$ the
indices are $S_{(ij)(ik)(jk)}$ with $i<j<k$, and if $i,j,k,l$ is a
permutation of $\{1,2,3,4\}$ we can use the shorthand notation
$S_l:=S_{(ij)(ik)(jk)}$, similarly $\tilde S_l:=\tilde
S_{(jk)(ik)(ij)}$. For $Q$'s we use $Q_i:=Q_{(1i)(jk)}$, where $j<k$
and $i,j,k$ is a permutation of $\{2,3,4\}$ (recall also Eqn.\
(\ref{E:QQ})). After multiplying each line from left and right by
suitable $\tilde S^{-1}$ and $Q$ to eliminate all $\tilde S$'s and
$Q^{-1}$'s,  exchanging left and right hand sides in
equations 3, 5 and 8, and writing the whole set in a different order
yields
\be
\begin{array}{ccc@{\hspace{1pt}}c@{\hspace{1pt}}c@{\hspace{1pt}}c
@{\hspace{1pt}}c@{\hspace{1pt}}ccc@{\hspace{1pt}}c@{\hspace{1pt}}c
@{\hspace{1pt}}c@{\hspace{1pt}}c@{\hspace{1pt}}c}
1,6\quad&&Q_4&S_4&S_3&Q_3Q_2&S_2&S_1&=&S_1&S_2&Q_2Q_3&S_3&S_4&Q_4,\\
7,5r\quad&&Q_4&\tilde S_4^{-1}&S_3&Q_3Q_2&S_2&S_1&=&
S_1&S_2&Q_2Q_3&S_3&\tilde S_4^{-1}&Q_4,\\
4,3r\quad&&Q_4&\tilde S_4^{-1}&\tilde S_3^{-1}&Q_3Q_2&S_2&S_1&=&
S_1&S_2&Q_2Q_3&\tilde S_3^{-1}&\tilde S_4^{-1}&Q_4,\\
2,8r\quad&&Q_4&\tilde S_4^{-1}&\tilde S_3^{-1}&Q_3Q_2&\tilde S_2^{-1}&S_1&=&
S_1&\tilde S_2^{-1}&Q_2Q_3&\tilde S_3^{-1}&\tilde S_4^{-1}&Q_4,\\
1r,6r\quad&&Q_4&\tilde S_4^{-1}&\tilde S_3^{-1}&Q_3Q_2&\tilde S_2^{-1}
&\tilde S_1^{-1}
&=&\tilde S_1^{-1}&\tilde S_2^{-1}&Q_2Q_3&\tilde S_3^{-1}&
\tilde S_4^{-1}&Q_4,\\
5,7r\quad&&Q_4&S_4&\tilde S_3^{-1}&Q_3Q_2&\tilde S_2^{-1}&
\tilde S_1^{-1}&=&
\tilde S_1^{-1}&\tilde S_2^{-1}&Q_2Q_3&\tilde S_3^{-1}&S_4&Q_4,\\
3,4r\quad&&Q_4&S_4&S_3&Q_3Q_2&\tilde S_2^{-1}&\tilde S_1^{-1}&=&
\tilde S_1^{-1}&\tilde S_2^{-1}&Q_2Q_3&S_3&S_4&Q_4,\\
8,2r\quad&&Q_4&S_4&S_3&Q_3Q_2&S_2&\tilde S_1^{-1}&=&
\tilde S_1^{-1}&S_2&Q_2Q_3&S_3&S_4&Q_4.
\eearn\label{E:13}\ee
This is the final form of the equations, when considered from the
algebraic point of view.  Clearly, if the unitarity condition
(\ref{E:Uni}) holds we have $\tilde S_i^{-1}=S_i$ and all equations
are identical.

Equations of exactly the same form but with $Q=\delta\delta$ were
presented in \cite{KV2}, where they were derived in quite a different
context: The study in \cite{KV2} was based on bicategories and
the equations were presented as a theorem stating that certain
bicategories will satisfy this set.

\subsection{Geometric interpretation}
Above we have given an algebraic derivation, but the tetrahedron
equations can also be derived by other approaches, for example by the
geometric approach of straight string scattering, which we will
consider next. In this approach the unitarity condition does not
arise, and we still have the ordering problem discussed in the
introduction. (In the following we will drop the $Q$-matrices.)

In order to have a geometric interpretation of (\ref{E:l1}) we
renumber the indices so that in each equation the indices of an
$S$-matrix contain the same set of numbers as in Case 1. The required
renumberings are: 1,6 none, 2,8r none, 3,4r (1324) and exchange of
left and right hand sides, 4,3r (132), 5,7r (1234) and exchange, 1r,6r
none, 7,5r (123), 8,2r (14) and exchange, this yields
\bearn
1,6&&\ST{12}{13}{23}\ST{12}{14}{24}\ST{13}{14}{34}\ST{23}{24}{34}
=\ST{23}{24}{34}\ST{13}{14}{34}\ST{12}{14}{24}\ST{12}{13}{23}\\
2,8r&&\ST{12}{13}{23}\tilde \ST{12}{42}{41}\tilde \ST{13}{43}{41}
\tilde \ST{23}{43}{42}
=\tilde \ST{23}{43}{42}\tilde \ST{13}{43}{41}\tilde
\ST{12}{42}{41}\ST{12}{13}{23}\\ 3,4r&&\tilde
\ST{12}{32}{31}\tilde \ST{12}{42}{41}\ST{34}{31}{41}\ST{34}{32}{42}
=\ST{34}{32}{42}\ST{34}{31}{41}\tilde \ST{12}{42}{41}\tilde
\ST{12}{32}{31}\\ 4,3r&&\ST{31}{32}{12}\ST{41}{42}{12}\tilde
\ST{41}{31}{34}\tilde \ST{42}{32}{34} =\tilde \ST{42}{32}{34}\tilde
\ST{41}{31}{34}\ST{41}{42}{12}\ST{31}{32}{12}\\ 5,7r&&\tilde
\ST{31}{21}{23}\tilde \ST{41}{21}{24}\tilde \ST{41}{31}{34}
\ST{23}{24}{34}
=\ST{23}{24}{34}\tilde \ST{41}{31}{34}\tilde \ST{41}{21}{24}\tilde
\ST{31}{21}{23}\\ 1r,6r&&\tilde \ST{23}{13}{12}\tilde
\ST{24}{14}{12}
\tilde \ST{34}{14}{13}\tilde \ST{34}{24}{23}=\tilde \ST{34}{24}{23}
\tilde \ST{34}{14}{13}\tilde \ST{24}{14}{12}\tilde \ST{23}{13}{12}\\
7,5r&&\ST{23}{21}{31}\ST{24}{21}{41}\ST{34}{31}{41}\tilde \ST{34}{24}{23}
=\tilde \ST{34}{24}{23}\ST{34}{31}{41}\ST{24}{21}{41}\ST{23}{21}{31}\\
8,2r&&\tilde \ST{23}{13}{12}\ST{41}{42}{12}\ST{41}{43}{13}\ST{42}{43}{23}
=\ST{42}{43}{23}\ST{41}{43}{13}\ST{41}{42}{12}\tilde \ST{23}{13}{12}
\eearn\label{E:l4}
\ee
These orderings were derived algebraically but one can give a geometric
rule that produces the same:

\noindent{\bf Geometric rule for label ordering:} 
Draw a line on the plane, outside the region of string intersections.
When the line moves, without changing its direction, it will sweep
across the intersection region.  {\em For each scattering matrix write
the indices of the corresponding triangle corners in the order the
line hits them}. If the order is counterclockwise, use $\tilde S$.

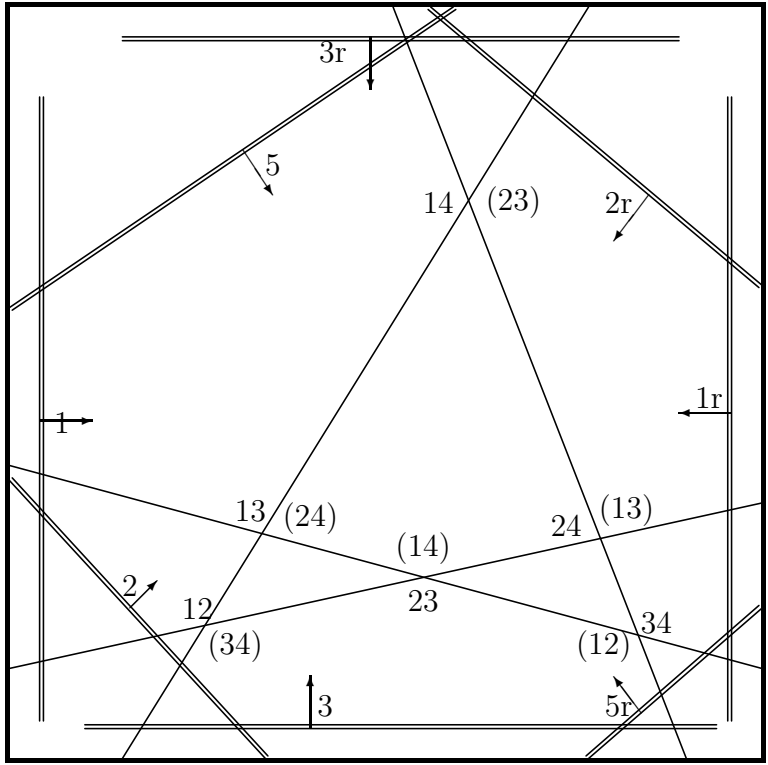
\begin{figure}[t]
\unitlength=1mm
\special{em:linewidth 0.6pt}
\linethickness{1.5pt}
\begin{picture}(110,100)(-10,10)
\put(10,10){\framebox(100,100)[cc]{}}
\emline{14}{15}{9}{14}{98}{10}
\emline{14.5}{15}{9}{14.5}{98}{10}
\emline{10}{70}{11}{69}{110}{12}
\emline{10.35}{69.65}{11}{69.35}{109.65}{12} 
\emline{25}{106}{13}{99}{106}{14}
\emline{25}{105.5}{13}{99}{105.5}{14}
\emline{110}{73}{15}{66}{110}{16}
\emline{109.65}{72.65}{15}{65.65}{109.65}{16}
\emline{106}{98}{17}{106}{15}{18}
\emline{105.5}{98}{17}{105.5}{15}{18}
\emline{87}{10}{19}{110}{30}{20}
\emline{86.65}{10.35}{19}{109.65}{30.35}{20}
\emline{20}{14}{21}{104}{14}{22}
\emline{20}{14.5}{21}{104}{14.5}{22}
\emline{44}{10}{23}{10}{47}{24}
\emline{44.35}{10.35}{23}{10.35}{47.35}{24}
\linethickness{0.6pt}
\emline{25}{10}{1}{87}{110}{2}
\emline{100}{10}{3}{61}{110}{4}
\emline{110}{22}{5}{10}{49}{6}
\emline{10}{22}{7}{110}{44}{8}
\put(14,55){\vector(1,0){7}}
\put(41,91){\vector(2,-3){4}}
\put(58,106){\vector(0,-1){7}}
\put(95,85){\vector(-3,-4){4.67}}
\put(106,56){\vector(-1,0){7}}
\put(94,16){\vector(-3,4){3.67}}
\put(50,14){\vector(0,1){7}}
\put(26,30){\vector(1,1){3.67}}
\put(35,30){\makebox(0,0)[cc]{12}}
\put(40,25){\makebox(0,0)[cc]{(34)}}
\put(42,43){\makebox(0,0)[cc]{13}}
\put(50,42){\makebox(0,0)[cc]{(24)}}
\put(67,84){\makebox(0,0)[cc]{14}}
\put(77,84){\makebox(0,0)[cc]{(23)}}
\put(65,31){\makebox(0,0)[cc]{23}}
\put(65,38){\makebox(0,0)[cc]{(14)}}
\put(84,41){\makebox(0,0)[cc]{24}}
\put(92,43){\makebox(0,0)[cc]{(13)}}
\put(96,28){\makebox(0,0)[cc]{34}}
\put(89,25){\makebox(0,0)[cc]{(12)}}
\put(17,56){\makebox(0,0)[ct]{1}}
\put(45,89){\makebox(0,0)[cc]{5}}
\put(53,104){\makebox(0,0)[cc]{3r}}
\put(91,84){\makebox(0,0)[cc]{2r}}
\put(103,58){\makebox(0,0)[cc]{1r}}
\put(91,17){\makebox(0,0)[cc]{5r}}
\put(52,17){\makebox(0,0)[cc]{3}}
\put(26,33){\makebox(0,0)[cc]{2}}
\end{picture}
\caption{How the possible algebraic orderings can be derived
geometrically.  The vertices are ordered according to the order in
which the moving (double) lines hit them. The complementary numbering
given in parentheses is used with the Bruhat ordering of
Sec. \protect\ref{bsec}}
\label{F:3}
\end{figure}
In Figure 3 we have redrawn Figure \ref{F:1} with 8 approaching lines.
These sweeping lines give exactly the orderings that were obtained by the
algebraic method after relabeling (\ref{E:l4}).

\subsection{Connection with Bruhat order B(4,2)}
\label{bsec}
The above concrete geometrical approach can be made more precise using
the notion of higher Bruhat orders, introduced by Manin and Shechtman
in \cite{MS,MS2}. This provides the proper algebraic setting for the
description of the general $d$-simplex equations, viewed as
higher-order intertwining or braiding objects. The setting is that of
moves of hyper-planes embedded in a higher-dimensional space.
Realizations can be constructed in terms of generators of the
fundamental group of configuration spaces formed by the complements of
such configurations \cite{Orlik}. An explicit realization was
constructed by Lawrence in \cite{Lawr2}.

To make our account self-contained, we will briefly describe the
Manin-Shechtman construction, cf. \cite{MS,MS2}.  (An alternative
description was given recently in \cite{Z}.)  First, for any pair of
integers, ~$n,k$, with $n\geq k\geq 1$, they introduce the set of
$k$-element subsets $C(n,k)$ of the set ~$\underline{n}=\{ 1,2,\dots
,n\}$, whose elements will be denoted by ~$( i_1 i_2\dots i_k)$~ in
increasing order, ~$i_1<i_2<\dots <i_k$~. For any given element $c\in
C(n,k)$, we denote by $\hat{c}_j$ the element of $C(n,k-1)$ obtained
by removing the $j$th element $i_j$, ($1\leq j\leq k$), from the tuple
$c$.

Next, we consider the set of total orders on the set $C(n,k)$.  For
this purpose we need to select from $C(n,k)$ only those orderings that
descend from $C(n,k+1)$, i.e. the elements $\hat{d}_j\in C(n,k)$,
($j=1,\dots,k+1$), coming from {\em each }$d\in C(n,k+1)$ by applying
the $\hat{}$-operation described above. They are ordered in either
ascending or in descending order,
\be
 \hat{d}_1<\hat{d}_2< \dots < \hat{d}_{k+1}\   \mbox{ or }
 \hat{d}_1>\hat{d}_2> \dots > \hat{d}_{k+1}\  . 
\label{eq:ddd} 
\ee
The set of all such total orderings is called $A(n,k)$, and its 
elements can be written as chains $a=c_1c_2\dots c_n$, $c_i\in 
C(n,k)$, with ~$c_1< c_2<\dots c_n$~ in the given ordering by $a$ . 

\paragraph{Example:} Consider the set $A(4,2)$ which are constructed 
according to the scheme above. First we need the sets 
$C(4,2)$ having six elements: 
\[
C(4,2)= \{ (12),(13),(14),(23),(24),(34) \}=:\{c_1,c_2,c_3,c_4,c_5,c_6\},
\]
and $C(4,3)$ having four elements: 
\[ C(4,3)= \{ (123),(124),(134),(234) \}\  .\] 
{}From the latter set we can construct the elements 
$\hat{d}_j$ for each $d\in C(4,3)$, leading to the following list of 
conditions on the orderings: 
\begin{eqnarray} 
(23)<(13)< (12) \   \ & \mbox{ or }& \    \ (23)>(13)> (12)\  , 
\nonumber \\  
(24)<(14)< (12) \   \ & \mbox{ or }& \    \ (24)>(14)> (12)\  , 
\nonumber \\  
(34)<(14)< (13) \   \ & \mbox{ or }& \    \ (34)>(14)> (13)\  , 
\nonumber \\  
(34)<(24)< (23) \   \ & \mbox{ or }& \    \ (34)>(24)> (23)\  . 
\nonumber   
\end{eqnarray} 

A geometric picture is very useful to find out which orderings on
$C(4,2)$ (i.e., which combinations of the above list) are actually
allowed. It turns out that the allowed orderings are exactly those
that can be obtained from figure \ref{F:3}, by looking at it from
various directions.  At this point we should keep also those orderings that
differ only by an exchange of elements not directly connected, e.g, by
a small tilt on direction 3 we can have $c_4c_2c_5\dots$ and
$c_4c_5c_2\dots$ In this way we obtain the set $A(4,2)$, containing
e.g.
\[
 A(4,2)= \{ c_1c_2c_3c_4c_5c_6\,,\, 
c_4c_2c_5c_3c_1c_6\,,\,c_6c_1c_3c_5c_2c_4\,,\, \dots \}\  . 
\] 

To obtain the Bruhat orders $B(n,k)$, we need to consider the set
$A(n,k)$ up to an inversion, i.e. selecting only one of each
possibility in (\ref{eq:ddd}). So, the set $Inv(a)$ of inversions on
an element $a\in A(n,k)$ is a subset $d\in C(n,k+1)$ such that
~$\hat{d}_1<\hat{d}_2<\dots <\hat{d}_{k+1}$~.  Furthermore, we need to
introduce an equivalence under adjacency.  Two elements in $A(n,k)$
will be called equivalent, ~$a\sim a^\prime$~, provided $a^\prime$ is
obtained from $a=c_1c_2\dots c_n$ by the permutation of two adjacent
subsets $c_j$ and $c_{j+1}$, containing in the union at least k+2
elements. The Bruhat orders are then the equivalence classes in
$A(n,k)$ under this equivalence relation, i.e. they are contained in
the set ~$B(n,k)=A(n,k)/\sim$~.
 
Let us now see what this amounts to in the case of $A(4,2)$.  The
adjacent elements in $C(4,2)$ are exactly the ones that have no
entries in common. In this way they correspond to the orderings up to
interchanging the subsets $(12)$ and $(34)$, $(13)$ and $(24)$,
and $(14)$ and $(23)$. Thus, with the above identification of orders,
we get for $B(n,k)$:
\begin{eqnarray} 
B(4,2) &=& \{ [c_1c_2c_3c_4c_5c_6]\,,\, [c_6c_5c_4c_3c_2c_1]
\,,\,[c_4c_2c_1c_3c_5c_6]\,,\, [c_1c_2c_3c_6c_5c_4]  \nonumber \\ 
&&\phantom{ \{ }[c_4c_5c_6c_3c_2c_1]\,,\,[c_6c_5c_3c_1c_2c_4]\,,\,
[c_6c_1c_3c_2c_5c_4]\,,\,[c_4c_5c_2c_3c_1c_6] \}\ ,
\end{eqnarray} 
corresponding to directions $1r,\,1,\,2r,\,5r,\,5,\,2,\,3,\,3r$ in
Fig.~\ref{F:3}.

It is easily noted that this partial ordering when imposed on 
$B(4,2)$ corresponds exactly to the  configurations in the 
obstruction derivation of the eight tetrahedron equations. 
At this point we can note the connection with the work on 
braided monoidal 2-categories, cf. \cite{VK}, that also leads 
to the set of eight tetrahedron equations, albeit from quite a 
different point of view.

\section{Conclusions}
The statement of this paper is the following: what  is usually
referred to as the tetrahedron equation is actually one of a system
of eight coupled equations that can be derived systematically from the
collection of all consistency conditions that arise from the underlying
set of trilinear equations (\ref{E:rev1}) and (\ref{E:rev2}).
We have shown also that the various
classes of starting configurations that lead to these  different
equations are labeled by a new algebraic object, which is the higher
Bruhat order $B(4,2)$ introduced by Manin and Shechtman. It is
obvious that these considerations can, in principle, be extended to
any dimension, i.e. to obtain systems of $D$-simplex equations for any
$D=2,3,\dots$~. 

Of course, our derivation comes down to the same type of combinatorics
that is behind the description in terms of 2-category theory,
\cite{VK}, but ours is closer to the physical interpretation.
Furthermore, we hope that the obstruction derivation might eventually
lead to the derivation of solutions to the system of eight equations
(in the cases that the set cannot be reduced to one single equation)
along the same lines as the derivation of solutions in the papers by
Korepanov, \cite{Kor1,Kor2}. We have investigated the known solutions
of the tetrahedron equations in \cite{Kor1,Labels}, but unfortunately
these all fall into the class of unitary solutions for which the
system collapses (that is, the second equation of (\ref{E:13}) did not
have other solutions than those with $\tilde S_4^{-1}=S_4$). However,
it cannot be ruled out that nontrivial solutions of the full
non-degenerate system (\ref{E:13}) exist even though it might not be so
easy to find such solutions.

\end{document}